---------------------------------------------------------------------------
############################## .TEX FILE ###################################
---------------------------------------------------------------------------

\documentstyle[12pt]{article}    
\begin{document}

\begin{center}
{\large \bf  A New Definition of Hypercomplex Analyticity}
\end{center}
\begin{center}
      {\em Stefano De Leo}\footnote{Supported in part by EEC contract, 
                              grant no.~ERBCHRX CT93-413}$^{\; (a,b)}$
        and
       {\em Pietro Rotelli}$^{\; (a)}$\\
        $^{(a)}$ Dipartimento di Fisica, 
                Universit\`a degli Studi Lecce\\ INFN, 
                Sezione di Lecce\\
                via Arnesano, CP 193, I-73100 Lecce, Italy\\
         $^{(b)}$ Dipartimento di Metodi e Modelli Matematici per le Scienze
                Applicate,\\
                via Belzoni 7, I-35131 Padova, Italy\\
January, 1997
\end{center}

\begin{center}
{\bf Abstract}
\end{center}
Complex analyticity is generalized to hypercomplex functions, 
quaternion or octonion, in such a manner that it includes
the standard complex definition and does not reduce analytic 
functions to a trivial class. A brief comparison with other 
definitions is presented.\\
{\tt 1991 Mathematics Subject Classifications}: 11R52, 30G35, 46S20.\\
{\tt Keywords}: Hypercomplex Numbers, Analytic Functions, 
Differential Operators.

	In this paper we shall propose a new definition for the 
generalization of the condition for complex analyticity for complex 
functions to quaternionic or octonionic functions of the corresponding 
variables. However, in the following unless explicitly mentioned, we 
will limit the discussion for clarity to quaternionic functions of 
quaternionic variables. 

	Since the discovery of quaternions by Hamilton~\cite{H43} in the 
last century a recurring question has been the best way to extend complex
analyticity to quaternionic functions of quaternionic variables. The most 
immediate idea is to extend the concept of differentiability from the
two-dimensional complex plain to the four-dimensional quaternionic space.
This can easily be done and involves the imposition of three quaternionic 
equations. If one defines a general quaternion by 
\begin{equation}
	q = x_{0} + i x_{1} + j x_{2} + k x_{3}~,
\end{equation}
where, $x_{0,1,2,3} \in {\cal R}$ and $i$, $j$ and $k$ are the noncommuting 
imaginary units such that  
\[ i^{2}=j^{2}=k^{2}=-1 ~~~~~~ \mbox{and} ~~~~~~ [i, \, j]=2k ~~~~ 
\mbox{(cyclic)}~,
\]
then the condition for differentiability of the quaternionic function 
$f(q)$ yields the three equations:
\begin{equation}
\label{three}
\partial_{x_{0}}f(q)=-i \partial_{x_{1}} f(q)=-j \partial_{x_{2}}f(q)=
-k\partial_{x_{3}}f(q)~.
\end{equation} 
The result of these equations is the restriction of the class of 
differentiable functions to merely linear functions in $q$. 
\begin{equation}
f(q)=c_1 + qc_2 ~ , 
\end{equation}
where $c_{1,2}$ are arbitrary quaternionic coefficients. 
Actually, the position of the constant depends upon whether the inverse 
of the increment in the derivative is placed upon the left or the right of 
the increment of the function $f(q)$. The above form assumes the left choice. 
If both options are imposed simultaneously then the holomorphic class of 
functions reduces to linear polinomials of $q$ with $c_{2}$ real.   
In any case this generalization is so restrictive that it loses all 
practical interest. We also note, that it excludes (except for the trivial
constant functions) the class of analytic/holomorphic functions of standard 
complex analysis, characterized by Taylor series in the corresponding 
complex variable. 

	A more sophisticated attempt was made in the 1930's by Fueter
and collaborators~\cite{F35}, who defined analyticity by means 
of a single quaternionic partial differential equation which includes 
the standard complex analyticity equation for complex functions of the 
corresponding single complex variable. He showed that this definition led 
to close analogues of Cauchy's theorem, Cauchy's integral formula, and the 
Laurent expansion. A complete bibliography of Fueter's work is contained 
in~\cite{HAE} and a simple account of the elementary parts of the theory 
has been given by Deavours~\cite{DEA}. A further extension of Fueter 
involves a third-order differential equation~\cite{GUR} which we mention 
only briefly in this work. However, if either of these is the best choice 
is still to be demonstrated, and from time to time variations on the 
theme appear in the literature~\cite{GUR,SUD,GT}. The situation is 
enriched by the non commutative nature of hypercomplex numbers which 
permits the definition of so-called left/right derivative operators 
according to the position of the imaginary units with respect to the 
function. Even combinations or admixtures of these alternative 
derivatives may be used. This ``complication'' is significant in certain 
other physical and mathematical applications such as in quaternionic group 
theory~\cite{GUT} and in the representation of Lorentz 
transformations~\cite{REL}.

	 The proposal in this work is much simpler. It is to define a 
``local'' derivative operator that depends upon the four-dimensional 
point at which the derivative is to be made. Each non-real quaternion 
point together with the real axis defines a unique complex plain and 
it is the complex variable of this plain that we use to define analyticity, 
in complete analogy to complex analysis. As a consequence, the class 
of analytic functions are generalized to include all polynomial functions 
of a single quaternionic variable with (right acting) quaternionic 
coefficients.  This approach has the non trivial virtue of being 
directly generalizable to octonionic functions of octonionic variables.    

	The complex derivative operator, $\partial_z$, is defined 
for $z \,=\, x+iy $, by
\begin{equation}
\partial_z = \frac{1}{2} \; (\partial_x - i\partial_y)~.
\end{equation}
Such that its action upon a monomial of $z$ is simply,
\begin{equation}
	\partial_z z^n = n z^{n-1}
\end{equation}
while it gives zero if applied to a polynomial of $\bar{z}=x-iy$.
Similarly the derivative operator for $\bar{z}$ can be defined by the 
operator
\begin{equation}
	\partial_{\bar{z}}  = \frac{1}{2} \; (\partial_x +i\partial_y)~.
\end{equation}
The conditions for a regular function $f(z,\bar {z})$, the Cauchy-Riemann 
equations, can then be expressed by
\begin{equation}
	\partial_{\bar{z}}f(z,\bar {z})=0
\end{equation}
The well known solutions of this equation are polynomial functions 
of $z$ without any dependence upon $\bar{z}$. At the level of complex 
numbers this definition of analyticity coincides with the existence of 
a unique derivative.

	The natural generalization of the above complex derivative 
operator to a quaternionic derivative, acting from the left, is the 
operator,
\begin{equation}
\partial_q  =  \frac{1}{2} \; \left( 
 \partial_{x{_0}}  - \frac{i\partial_{x{_1}} + j\partial_{x{_2}} 
+ k\partial_{x{_3}}}{3} \right)   ~.  
\end{equation}
We have chosen the normalization factors in $\partial_q$  so that it gives 
the expected result when applied to linear 
functions of the variable $q$, i.e.
\begin{equation}
\partial_q (c_1+q c_2) = c_2~,
\end{equation}
with $c_1$ and $c_2$  quaternionic constants as before. This operator also 
annihilates 
\[ \bar{q}= x_0 - ix_1 - jx_2 - kx_3 \]
However, it does not act in a simple way upon higher powers of $q$ or 
$\bar {q}$. More importantly, it does not reduce to the complex derivative 
operator when applied to functions independent of, say, $x_{2}$ and 
$x_{3}$.

	One of the possible alternatives, at this point, is to define 
$\partial_{\bar{q}}$  so that it yields the desired complex limit.
This limit corresponds to functions of only the real and one imaginary 
variable and involving only one imaginary unit (a limitation upon the
constants).  The proposal of Fueter is to define this operator by,
\begin{equation}
\label{feuter}
\partial_{\bar {q}} = \partial_{x_{0}} + i\partial{x_{1}} + 
j\partial_{x_{2}} + k\partial_{x_{3}} 
\end{equation}
This yields the so called Cauchy-Riemann-Fueter (CRF) equation,
\begin{equation}
	\partial_{\bar{q}}f = 0~, 
\end{equation}
where $\bar{q}$ is the conjugate of $q$ and $\partial_{\bar{q}}$ the 
conjugate of $\partial_q$. Not even linear functions in $q$ satisfy 
this equation. We seem to have reduced analytic functions to mere 
constants. But one quaternionic equation is surely less restrictive 
than the three of holomorphy first quoted, Eq.~(\ref{three}). 
In fact, the CRF equation has many (but not $q$ polinomial) 
solutions particularly in what we shall call the {\em canonical complex
variable  limit} (see below). In complex analysis, there is no analogy for the 
appearance of 
new solutions when one reduces the number of variables. There, if a 
function depends upon e.g., only the real variable, the analyticity 
conditions (Cauchy-Riemann equations) reduce the function to a constant 
which is already included in the general class of polinomial functions 
in $z$.

     The canonical complex limit refers to polinomial functions of 
complex variables involving {\em only}  one of the {\em basic} imaginary units, 
i.e. with as variable $z$ either $x_0 + ix_1$ or $x_0 + jx_2$ or 
$x_0 + kx_3$. These are regular in the above Fueter sense. 
Nevertheless, since we can multiply each term in the Taylor expansion 
of the analytic function {\em from the right} by independent 
quaternionic coefficients and still maintain analyticity, $f(z)$ 
need not lie in the same complex plain as its variable $z$, nor 
indeed be restricted to lie in a particular plain in quaternionic space. 
Nevertheless, with the above restriction 
to {\em canonical complex variables}, the standard complex analytic functions 
with complex coefficients are indeed CRF quaternionic analytic.  On the 
other hand, a complex ($1, \, \iota $) polinomial of a complex 
$\zeta =x_0 + \iota x$, 
where $\iota$ is a linear (not necessarily canonical) normalized combination of 
$ix_1$, $jx_2$ and $kx_3$, is {\em not} in general  a regular function, 
in the Fueter sense. 
We also observe that the above mentioned functions by no means 
exhaust the class of quaternionic analytic functions, as can easily be seen 
by considering the set of functions independent of only one of the real 
variables such as $x_3$. 
This definition of regular functions by Fueter 
permits a quaternionic version of Cauchy's theorem and of Cauchy's integral 
formula.
	We would like know to define an analyticity condition which both 
includes arbitrary complex functions of corresponding complex 
variables projected from $q$, and even more important, extends 
the class of regular functions to polinomial functions in $q$. 
To the best of our knowledge, the method we describe below has not been 
proposed previously. The {\em trick} is obvious once seen. It however requires 
the extension of the concept of a derivative operator from those with 
constant (even quaternionic) coefficients, to one with variable 
coefficients dependent upon the point of application. In physical terms 
we would say (in analogy with gauge transformations) that the 
derivative operator passes from a {\em global} form to a {\em local} form. 

   As already mentioned in the introduction, a point $q \not \in {\cal R}$ 
together with the real axis $q=x_0$ specifies a unique 
complex plain with imaginary axis given by the unit
\begin{equation}
\iota  = (i x_1 + j x_2 + kx_3)/ |\vec{x}|~,
\end{equation}
where $|\vec{x}|^2 = x_1^2 + x_2^2 + x_3^2$. With $\zeta$ already defined 
as $x_0 + \iota x$, and $x=|\vec x|$, the analyticity condition we propose 
reads,
\begin{equation}
 \partial_{\bar{\zeta}} f(\zeta) =
\frac{\partial_{x_0} + \iota \partial_x}{2} \, f(q)= 0~,
\end{equation}
or explicitly, ignoring the constant factor of $1/2$,
\begin{equation}
\left[ \partial_{x_0} + \frac{ix_1+jx_2+kx_3}{|\vec{x}|^2} 
\, (x_1 \partial_{x_1}+x_2 \partial_{x_2}+x_3 \partial_{x_3}) \right] 
f(q) = 0~.
\end{equation}
This formally reproduces to the Cauchy-Riemann equations in the 
corresponding complex plain, which include, as special cases, the three 
canonical plains. 
However, it involves a quaternionic derivative and includes amongst the 
class of regular functions {\em all} polynomial functions of $q$ with 
right-acting quaternionic 
coefficients. That is, since 
\begin{equation} 
 [~\partial_{\bar{\zeta}}, \; q~] =  0~,
\end{equation}
a class of solutions of the analyticity conditions is the Taylor series
\begin{equation}
f(q)=\sum_{n\geq {0}} q^{n}c_{n}
\end{equation}
where the $c_{n}$ are arbitrary quaternionic constants. 
This is the same class of functions obtained by Fueter by means of a 
{\em third order} analyticity condition. A more detailed comparison 
of our approach with the ever green work of Fueter will be presented 
elsewhere~\cite{DR}. It should be obvious that a similar definition 
of a regular function for octonions reproduces a similar result. 
The non associative nature of octonions plays no role in this since 
our analyticity condition is of first order.
 
	In conclusion we have described, from what we hope is an original 
viewpoint, the quest for an appropriate quaternionic analyticity condition.
We have introduced a local derivative operator and used it in a first order
analyticity condition which includes and generalizes to hypercomplex 
functions the equations of Cauchy-Riemann. The result is a rich and we 
believe significant class of analytic functions obtained previously 
but in a much more laborious manner by Fueter. 

The authors of this letter wish to dedicate this work to the memory of 
{\bf Abdus Salam} who was professor of one of us (PR) at Imperial 
college during the exciting years of the late 1960's. He promoted in 
his students a passion for both Physics and Mathematics and we hope 
that some of this passion can be in turn passed on to even younger 
generations.

%%%%%%%%%%%%%%%%%%%%%%%%%%%%%%%%%%%%%%%%%%%%%%%%%%%%%%%%%%%%%%%%%%%%%%%%%%%%%%%
%                                REFERENCES
%%%%%%%%%%%%%%%%%%%%%%%%%%%%%%%%%%%%%%%%%%%%%%%%%%%%%%%%%%%%%%%%%%%%%%%%%%%%%%%

%%%%%%%%%%%%%%%%%%%%%%%%%%%%%%%%%%%%%%%%%%%%%%%%%%%%%%%%%%%%%%%%%%%%%%%%%%%%%%
\end{document}